\documentclass{amsart}

\usepackage{amssymb}
\usepackage{eepic}
\usepackage{color}
\usepackage{mathtools,xparse}
\usepackage[linktocpage=true]{hyperref}
\usepackage[utf8x]{inputenc}

\allowdisplaybreaks

\numberwithin{equation}{section}

\newtheorem{theorem}{Theorem}[section]
\newtheorem{lemma}[theorem]{Lemma}

\newtheorem{corollary}[theorem]{Corollary}
\newtheorem{remark}[theorem]{Remark}

\newtheorem{Ctheorem}{Theorem}

\newtheorem{Cremark}[Ctheorem]{Remark}
\newtheorem{Cconjecture}[Ctheorem]{Conjecture}

\newtheorem{TheoA}{Theorem A}

\newtheorem{TheoAAA}{Theorem B}

\newcommand{\Z}{\mathbf{Z}}

\newcommand{\R}{\mathbf{R}}
\newcommand{\C}{\mathbf{C}}

\newcommand{\B}{\mathcal{B}}
\newcommand{\SL}{S \hskip-1pt L_n(\R)}

\newcommand{\summ}{\sum\nolimits}

\def\G{\mathrm{G}}
\def\1{\mathbf{1}}

\def\Q{\mathcal{Q}}

\def\RR{\mathcal{R}}

\def\V{\mathrm{\mathcal{L}(G)}}

\newcommand{\dem}{\noindent {\bf Proof. }}
\newcommand{\ske}{\noindent {\bf Sketch of the proof. }}

\newcommand{\demA}{\noindent {\bf Proof of Theorem A. }}

\newcommand{\demAAA}{\noindent {\bf Proof of Theorem B. }}

\newcommand{\fin}{\hspace*{\fill} $\square$ \vskip0.2cm}

\def\mean{- \hskip-11pt \int}

\begin{document}

\null

\null

\begin{center}
{\huge A H\"ormander-Mikhlin theorem \\ [3pt] for higher rank simple Lie groups}

\vskip15pt

{\sc {Jos\'e M. Conde-Alonso, Adri\'an M. Gonz\'alez-P\'erez \\ Javier Parcet and Eduardo Tablate}}
\end{center}

\title[H\"ormander-Mikhlin theorem for simple Lie groups]{}


\maketitle

\null

\vskip-45pt

\null

\begin{center}
{\large {\bf Abstract}}
\end{center}

\vskip-25pt

\null

\begin{abstract}
We establish regularity conditions for $L_p$-boundedness of Fourier multipliers on the group von Neumann algebras of higher rank simple Lie groups. This provides a natural H\"ormander-Mikhlin criterion in terms of Lie derivatives of the symbol and a metric given by the adjoint representation. In line with Lafforgue/de la Salle's rigidity theorem, our condition imposes certain decay of the symbol at infinity. It refines and vastly generalizes a recent result by Parcet, Ricard and de la Salle for $\SL$. Our approach is partly based on a sharp local H\"ormander-Mikhlin theorem for arbitrary Lie groups, which follows in turn from recent estimates by the authors on singular nonToeplitz Schur multipliers. We generalize the latter to arbitrary locally compact groups and refine the cocycle-based approach to Fourier multipliers in group algebras by Junge, Mei and Parcet. A few related open problems are also discussed. 
\end{abstract}

\addtolength{\parskip}{+1ex}

\vskip20pt

\section*{\bf Introduction}

Regularity conditions for $L_p$-boundedness of Fourier multipliers are central in harmonic analysis, with profound applications in theoretical physics, differential geometry or partial differential equations. The H\"ormander-Mikhlin fundamental condition \cite{Ho,Mi} gives a criterion for $L_p$-boundedness of the Fourier multiplier $T_m$ associated to the symbol $m: \R^n \to \C$ $$\widehat{T_mf}(\xi) = m(\xi) \widehat{f}(\xi).$$ Namely, if $1 < p < \infty$ the following bound holds 
\begin{equation} \tag{HM} \label{Eq-HMClassic}
\big\| T_m \hskip-2pt: L_p(\R^n) \to L_p(\R^n) \big\| \, \lesssim \, \frac{p^2}{p-1} \hskip-2pt \sum_{|\gamma| \le [\frac{n}{2}] +1} \Big\| |\xi|^{|\gamma|} \big| \partial_\xi^\gamma m(\xi) \big| \Big\|_\infty.
\end{equation}
It imposes $m$ to be a bounded smooth function over $\R^n \setminus \{0\}$. Locally, it admits a singular behavior at $0$ with a mild control of derivatives around it up to order $[\frac{n}{2}] +1$. This singularity is linked to deep concepts in harmonic analysis and justifies the key role of the H\"ormander-Mikhlin theorem. The same derivatives decay asymptotically to $0$, at a polynomial rate dictated by the differentiation order. It is optimal in the sense that we may not consider less derivatives or larger upper bounds for them. A Sobolev type formulation admits fractional differentiability orders up to $\frac{n}{2} + \varepsilon$ for any $\varepsilon > 0$. Condition \eqref{Eq-HMClassic} up to order $\frac{n-1}{2}$ is necessary for radial $L_p$-multipliers and arbitrary $p < \infty$. A characterization of general Fourier $L_p$-multipliers is considered nowadays beyond the reach of Euclidean harmonic analysis methods.

The interest of Fourier multipliers over group von Neumann algebras was early recognized by Haagerup in his pioneering work on free groups \cite{H} and the research thereafter on semisimple lattices \cite{DCH,CH}, encoding deep geometric properties of these groups in terms of approximation properties. More recently, strong rigidity properties of higher rank lattices were found in the remarkable work of Lafforgue and de la Salle \cite{LdlS} studying $L_p$-approximations. This has motivated intense research on Fourier multiplier $L_p$-theory over group algebras \cite{GJP,JMP1,JMP2,dLdlS,MR,MRX,PRS,PRo} due to potential applications in geometric group theory and classification theory, particularly towards a possible solution of Connes' rigidity conjecture. A wide interpretation of tangent space for general topological groups was introduced in \cite{GJP,JMP1,JMP2} by means of finite-dimensional orthogonal cocycles $\beta: \G \to \R^n$. If $m: \G \to \C$ satisfies the identity $m = \widetilde{m} \circ \beta$, the main discovery was that a H\"ormander-Mikhlin theory in group von Neumann algebras is possible in terms of the $\beta$-lifted symbols $\widetilde{m}$. This unfortunately requires the action to be orthogonal and excludes infinite-dimensional cocycles. Moreover, it imposes auxiliary differential structures, which appear to be unnecessary or at least less natural for Lie groups. All of it justifies a great interest in H\"ormander-Mikhlin conditions on Lie group von Neumann algebras in terms of Lie differentiation and natural/intrinsic metrics.   

Let $\G$ be a unimodular Lie group and let $L_p(\V)$ be the noncommutative $L_p$ space over its group von Neumann algebra, equipped with its natural operator space structure \cite{P2}. Consider a Fourier multiplier $T_m$ associated to $m: \G \to \C$ and set $d_g^\gamma m(g)$ for its left-invariant Lie derivative of order $\gamma$ with respect to a fixed ONB in the Lie algebra. We investigate the inequality below for $1 < p < \infty$, some natural lengths $L: \G \to \R_+$ and certain constant $\Delta_\G$ depending on $\G$ \vskip-10pt 
\begin{equation} \tag{HM$_\G$} \label{Eq-HM}
\big\| T_m: L_p(\V) \to L_p(\V) \big\|_{\mathrm{cb}} \, \le \, C_p \sum_{|\gamma| \le \Delta_\G} \big\| L(g)^{|\gamma|} d_g^\gamma m(g) \big\|_\infty.
\end{equation} 
\noindent We focus on lengths $L(g) = \mathrm{dist}(g,e)$ coming from a natural metric on $\G$. Beyond \cite{PRS}, very little is known in this direction. Let us start by considering the length $L_\mathrm{R}$ which is inherited from the Riemannian metric on the Lie group $\G$. Is there a H\"ormander-Mikhlin (HM) criterion for general Lie group von Neumann algebras using Lie derivatives and the Riemannian metric? 

\begin{TheoA}[Local HM criterion] 
Let $\G$ be a $n$-dimensional unimodular Lie group equipped with its Riemannian metric $\rho$ and set $L_\mathrm{R}(g) =\rho(g,e)$. Let $1 < p < \infty$ and let $m: \G \to \C$ be a Fourier symbol supported by a relatively compact neighborhood of the identity $\Omega$. Then, the following inequality holds $$\big\| T_m: L_p(\V) \to L_p(\V) \big\|_{\mathrm{cb}} \, \le \, C_p(\Omega) \sum_{|\gamma| \le [\frac{n}{2}]+1} \big\| L_\mathrm{R}(g)^{|\gamma|} d_g^\gamma m(g) \big\|_\infty,$$ for left-invariant Lie derivatives $d_g^\gamma$ of order $\gamma = (j_1, j_2, \ldots, j_{|\gamma|})$ with $1 \le j_i \le n$.    
\end{TheoA} 

Theorem A is certainly optimal for Riemannian metrics. First, its local nature is necessary, as for simple Lie groups a global statement would get in conflict with Lafforgue/de la Salle's rigidity theorem \cite{LdlS}. Secondly, it follows from the classical HM theorem that the Mikhlin regularity order $\Delta_\G = [\frac12 \dim \G]+1$ is sharp for Lie derivatives and any locally Euclidean metric. Finally, for fixed $1 < p < \infty$ we have $C_p(\Omega) \approx p^2/(p-1)$ as $\Omega$ shrinks to the identity. This matches the behavior of the $L_p$-constant for the Hilbert transform, the archetype of HM multiplier. 

A highly technical argument in line with Theorem A was recently presented in \cite{PRS}, which led to nearly optimal regularity orders for special linear groups. It imposes though much higher regularity for other Lie groups with large Euclidean codimension. More precisely, it gives rise to $\Delta_\G \ge [N^2/2]+1$ for the minimal $N$ satisfying $\G \subset G\hskip-1pt L_N(\R)$, which is far from optimal when $\dim \G \ll N^2$. Our argument is much more efficient. The strategy consists in relating the problem with an equivalent formulation in terms of singular Schur multipliers. To do so we first use a local form of Fourier-Schur transference \cite{PRS}, which includes nonamenable groups as well. Then we lift the resulting Herz-Schur multiplier to a nonToeplitz multiplier in the Lie algebra via the exponential map, which leads to a Euclidean deformation of the original multiplier. This  opens a door to new $S_p$-boundedness criteria \cite{CGPT}. The goal then is to prove that our HM conditions may be deduced from the conditions for the lifted Schur multiplier on the Lie algebra. Locality is critical essentially in all steps of the strategy. Our second goal is to eliminate locality for simple Lie groups. 

Riemannian symmetric spaces are quotients $\G/\mathrm{K}$ of a real simple (noncompact and connected) Lie group $\G$ with a finite center and a maximal compact subgroup $\mathrm{K}$. Fourier multipliers for Riemannian symmetric spaces were first considered by Clerc and Stein \cite{CSt}. Stanton/Tomas in rank one \cite{ST} and Anker in higher ranks \cite{A} obtained optimal HM criteria in this context. Under less symmetric assumptions the dual problem on the whole simple Lie group (not quotients of it) was first considered in \cite{PRS} on the group algebra of $\G = S \hskip-1pt L_n(\R)$. The Introduction of that work provides a careful analysis of the connections of this subject with geometric group theory and operator algebra. Let us now consider any simple Lie group $\G$. As in \cite{PRS}, the natural length in this context $L_\G: \G \to \R_+$ is locally Euclidean around the identity and its asymptotic behavior is dictated by the adjoint representation $$L_\G(g) \approx \|\mathrm{Ad}_g\|^{\tau_\G} \quad \mbox{as} \quad g \to \infty$$ for $\tau_\G = \mathrm{d}_\G / [\frac12 (\dim \G + 1)]$ with $\mathrm{d}_\G$ from \cite{Mau}. We have $\tau_\G = \frac12$ when $\G = S \hskip-1pt L_n(\R)$.

\begin{TheoAAA}[HM for simple Lie groups]
Let $\G$ be a $n$-dimensional simple Lie group with $n \ge 2/\tau_\G$. Then, the following inequality holds for any Fourier symbol $m: \G \to \C$ and $1 < p < \infty$ $$\big\| T_m: L_p(\V) \to L_p(\V) \big\|_{\mathrm{cb}} \, \le \, C_p \sum_{|\gamma| \le [\frac{n}{2}]+1} \big\| L_\G(g)^{|\gamma|} d_g^\gamma m(g) \big\|_\infty.$$
\end{TheoAAA} 

The new difficulty in proving Theorem B relies on asymptotic behavior of the symbols. The local one follows from Theorem A, due to the Euclidean nature of the metric around the identity. Maucourant's constant $\mathrm{d}_\G$ gives the volume growth rate of $\mathrm{Ad}$-balls up to a logarithmic factor. We shall prove that the assumption $n \ge 2/\tau_\G$ ensures that our HM condition implies asymptotically $$|d_g^\gamma m(g)| \lesssim \|\mathrm{Ad}_g\|^{- \mathrm{d}_\G} \quad \mbox{for} \quad \mbox{$|\gamma| \le [\frac{n}{2}]+1$.}$$ A similar approach in \cite{PRS} imposed additional decay since the HM condition there was more restrictive, from which the asymptotic behavior could be deduced out of the local one by an elementary cut and paste argument. On the contrary, the optimal regularity order in Theorem B leads to a critical decay order, from where we need a more elaborated argument which (notably) fails below Mikhlin's critical regularity index. As noted in Remark \ref{Rem-Optimal}, this might indicate that there is no more room for improvement in the metric $L_\G$. In addition, we should also indicate that:    
\begin{itemize}
\item[i)] The parameter $\tau_\G$ is just relevant for the asymptotic behavior in Theorem B. In particular, the smaller $\tau_\G$ is the better is the metric and so is our HM condition. We already mentioned that $\tau_\G = \frac12$ for $\G = \SL$. What is perhaps more significant is that $\tau_\G \le 1$ for any simple Lie group $\G$ and Theorem B solves in great generality the problem initiated in \cite{PRS} with optimal regularity. By Fourier-Schur transference and restriction, similar bounds also hold for Herz-Schur multipliers on high rank lattices. 

\vskip5pt

\item[ii)] The opaque condition $n \ge 2/\tau_\G$ just means $$\mathrm{d}_\G \ge 2 \mbox{$[\frac12 (\dim \G + 1)]$}/\dim \G,$$ which holds for large classes of simple Lie groups, see \eqref{Eq-Mau}. It is worth mentioning though that such condition fails for $S \hskip-1pt L_2(\R)$, since $\mathrm{d}_{S\hskip-1pt L_2(\R)} = 1$ and $\dim S\hskip-1pt L_2(\R) =3$. The local behavior still works and improves \cite{PRS}, but we do not find fast enough decay of the symbol and its derivatives to follow our proof. This fast decay is necessary for higher rank Lie groups from the rigidity theorems in \cite{LdlS,PRS}. However, $S \hskip-1pt L_2(\R)$ is weakly amenable and we should expect to find Fourier multipliers with arbitrarily slow decay paces for it and other rank one simple Lie groups. In conclusion, new ideas seem to be necessary to understand this problem in depth.
\end{itemize}

The structure of the paper is the following. In Section \ref{Sect-HMS} we briefly review local Fourier-Schur transference from \cite{PRS} and our recent results for singular nonTopelitz Schur multipliers \cite{CGPT}. Sections \ref{Sect-Local} and \ref{Sect-Simple} include the proofs of Theorems A and B respectively. Section \ref{Sect-Groups} gives a generalization of the main result in \cite{CGPT} replacing Euclidean spaces by locally compact groups, and we shall apply it to generalize the main results in \cite{JMP1}. This includes nonorthogonal cocycles and nonunimodular groups, precise statements are given in the body of the paper. In the Appendix we consider a HM criterion for stratified Lie groups in terms of the subRiemannian metric, which seems very natural since its Euclidean analogue holds true. Then we discuss why the approach through Schur multipliers in the Lie algebra |included in a previous version of this paper| cannot hold. In view of the Levi decomposition for Lie algebras (in terms of solvable and semisimple subalgebras), HM criteria for stratified Lie groups and rank one simple Lie groups appear to be very interesting.  

\section{\bf Local transference and Schur multipliers} \label{Sect-HMS}

In this section we review local Fourier-Schur transference from \cite{PRS} and a recent H\"ormander-Mikhlin type inequality for nonToeplitz Schur multipliers \cite{CGPT}. Schur multipliers on $M_n(\C)$ are linear maps $S_M(A) = ( M(j,k) A_{jk} )_{jk}$ for some symbol $M: \{1,2,\ldots,n\} \times \{1,2,\ldots,n\} \to \C$. More general index sets correspond to operators $A$ on $L_2(\Omega,\mu)$ for some $\sigma$-finite measure space $(\Omega,\mu)$. Let $S_p(\Omega)$ be the Schatten $p$-class over $L_2(\Omega,\mu)$. Not every operator in the Schatten $p$-class admits a kernel or matrix representation, but this is the case in $S_p(\Omega) \cap S_2(\Omega)$. In particular we set $S_M(A) = (M(\omega_1,\omega_2) A_{\omega_1 \omega_2})$ for every operator $A$ in $S_2(\Omega)$ and we say that $S_M$ is completely $S_p$-bounded when it maps $S_p(\Omega) \cap S_2(\Omega)$ into $S_p(\Omega)$ and extends to a cb-map on the whole $S_p(\Omega)$. 

\subsection{Local transference}

Let $\G$ be a locally compact unimodular group with Haar measure $\mu$ and left regular representation $\lambda$. Its group von Neumann algebra $\V$ is defined as the weak-$*$ closure in $\mathcal{B}(L_2(\G))$ of $\mathrm{span}(\lambda(\G))$. We may approximate every element affiliated to $\V$ by operators $$f \, = \, \int_\G \widehat{f}(g) \lambda(g) \, d\mu(g)$$ for smooth enough $\widehat{f}$, see \cite[Appendix A]{JMP2}. If $e$ is the unit in $\G$, $\tau(f) = \widehat{f}(e)$ determines the Haar trace $\tau$. Given a symbol $m: \G \to \C$, its associated Fourier multiplier is the map $T_m: \lambda(g) \mapsto m(g) \lambda(g)$ which satisfies $$\widehat{T_m f}(g) = \tau(T_mf \lambda(g)^*) = m(g) \tau (f \lambda(g)^*) = m(g) \widehat{f}(g).$$ In other words, it intertwines pointwise multiplication with the Fourier transform. 

Let us now briefly review certain Fourier/Schur transference results. Given $\G$ a unimodular group and $m: \G \to \C$, consider its Fourier and Herz-Schur multipliers on the group and matrix algebras associated to $\G$. In other words, we formally have $$T_m(\lambda(g)) = m(g) \lambda(g), $$ \vskip-17pt $$S_m(A) = \big( m(gh^{-1})A_{gh} \big).$$ Here $\lambda: \G \to \mathcal{U}(L_2(\G))$ stands for the left regular representation. Finding accurate conditions on the symbol $m$ which ensure $L_p$-boundedness of these multipliers is a rather difficult problem. It is known from \cite{CS,NR} that both problems coincide for amenable groups. More precisely, if $S_p(\G)$ and $L_p(\V)$ stand for the natural noncommutative $L_p$ spaces on these algebras, we get
\begin{equation} \label{Eq-FS}
\hskip4pt \big\| S_m: S_p(\G) \to S_p(\G) \big\|_{\mathrm{cb}} \, = \, \big\| T_m: L_p(\V) \to L_p(\V) \big\|_{\mathrm{cb}} 
\end{equation}
for $\G$ amenable. Here, $\| \,\|_{\mathrm{cb}}$ denotes the completely bounded norm of these maps after endowing $S_p$ and  $L_p$ with their natural operator space structure \cite{P2}. Complete norms are essential in this isometry. The upper bound holds for nonamenable groups as well. The reverse inequality remains open. Nevertheless, a local form of it was recently proved in \cite[Theorem 1.4]{PRS} up to a constant measuring the amenability distortion. It will play an important role later in this paper.

\begin{theorem}[Local transference] \label{Thm-LocalTransf}
Let $\G$ be a locally compact unimodular group and consider a relatively compact neighborhood of the identity $\Omega$ and any open set $\Sigma$ in $\G$ containing the closure of $\Omega$. Let $m \hskip-2pt : \G \to \C$ be a bounded symbol supported in $\Omega$. Then, the following inequality holds for any $p \in 2\Z_+$ $$\big\| T_m \hskip-2pt : L_p(\V) \to L_p(\V) \big\|_{\mathrm{cb}} \, \le \, C_{\Omega,\Sigma,p} \, \big\| S_m \hskip-2pt : S_p(L_2(\Sigma)) \to S_p(L_2(\Sigma)) \big\|_{\mathrm{cb}}.$$ 
Moreover, we may put $C_{\Omega,\Sigma,p} \le 2$ for $\Omega$ small enough depending on the value of $p$. 
\end{theorem}

\begin{remark} \label{Rem-GralFSTransf}
\emph{It is worth mentioning that the upper bound in \eqref{Eq-FS} holds as well for nonunimodular groups using a more involved definition of Fourier multiplier \cite{CS}.}
\end{remark}

\begin{remark} \label{Rem-Linfty-FSTransf}
\emph{It is known from \cite{BF} that \eqref{Eq-FS} holds true when $p=\infty$ even for nonamenable groups. Moreover, cb-norms and classical norms coincide in this case.}
\end{remark}

\subsection{Singular Schur multipliers}

We now focus on Schur multipliers with index set $(\Omega,\mu)$ being the $n$-dimensional Euclidean space equipped with its Lebesgue measure. Other groups will be considered in Section \ref{Sect-Groups}. NonToeplitz multipliers $M(x,y) \neq m(x-y)$ are no longer related to Fourier multipliers as in \eqref{Eq-FS}. The following result from \cite{CGPT} gives a simple criterion for $S_p$-boundedness.

\begin{theorem}[HMS multipliers] \label{Thm-HMS} 
If $1 < p < \infty$, we have $$\big\| S_M \big\|_{\mathrm{cb}(S_p(\R^n))} \, \lesssim \, \frac{p^2}{p-1} \sum_{|\gamma| \le [\frac{n}{2}] +1} \Big\| |x-y|^{|\gamma|} \Big\{ \big| \partial_x^\gamma M(x,y) \big| + \big| \partial_y^\gamma M(x,y) \big| \Big\} \Big\|_\infty,$$ for every Schur symbol $M: \R^n \times \R^n \to \C$ in the class $C^{[n/2]+1}(\R^{2n} \setminus \{x=y\})$.
\end{theorem} 

Theorem \ref{Thm-HMS} introduces a new class of singular $S_p$-bounded multipliers which we shall refer to as H\"ormander-Mikhlin-Schur (HMS) multipliers. Observe that \eqref{Eq-FS} shows that the above condition for Toeplitz symbols $M(x,y) = m(x-y)$ reduces to the classical condition \eqref{Eq-HMClassic} on $m$. This will be a key tool below. 

\section{\bf The local theorem} \label{Sect-Local}

Let $\G$ be a unimodular Lie group, consider the left-invariant vector fields in $\G$ generated by an orthonormal basis $\mathrm{X}_1, \mathrm{X}_2, \ldots, \mathrm{X}_n$ of $\mathfrak{g}$. Then, the corresponding Lie derivatives are defined as  
$$\partial_{\mathrm{X}_j} m (g) \, = \, \frac{d}{ds}\Big|_{s=0} m \big( g \exp(s \mathrm{X}_j) \big)$$
and do not commute for $j \neq k$. This justifies to define the set of multi-indices $\gamma$ as ordered tuples $\gamma = (j_1, j_2, \ldots, j_k)$ with $1 \le j_i \le \dim \G$ and $|\gamma| = k \ge 0$, which correspond to the Lie differential operators $$d_g^\gamma m (g) \, = \, \partial_{\mathrm{X}_{j_1}} \partial_{\mathrm{X}_{j_2}} \cdots \, \partial_{\mathrm{X}_{j_{|\gamma|}}} m (g) \, = \, \Big( \prod_{1 \le k \le |\gamma|}^{\rightarrow} \partial_{\mathrm{X}_{j_k}} \Big) m(g).$$

\demA By interpolation and duality, it suffices to prove that the statement holds for any $p \in 2 \Z_+$. Let us fix such an even integer, let $\Omega, \Sigma \subset \G$ be relatively compact symmetric neighborhoods of the identity. According to local transference Theorem \ref{Thm-LocalTransf}, the following inequality holds when $\mathrm{supp}(m) \subset \Omega$ for sufficiently small $\Omega$
$$\big\| T_m: L_p(\V) \to L_p(\V) \big\|_{\mathrm{cb}} \, \lesssim \, \big\| S_m: S_p(\Sigma) \to S_p(\Sigma) \big\|_{\mathrm{cb}}.$$ Let $\Xi$ be an open set containing the closure of $\Sigma$ and let $\phi: \G \to \R_+$ be a smooth function  identically $1$ in $\Sigma$ and vanishing outside $\Xi$. The map $\exp: \mathfrak{g} \to \G$ is a local diffeomorphism at the origin of the Lie algebra $\mathfrak{g}$. Therefore, if we pick $\Xi$ small by taking $\Sigma$ and $\Omega$ small enough, this yields a diffeomorphism $\exp: \mathrm{U} \to \Xi$ over certain neighborhood $\mathrm{U}$ of the origin of $\mathfrak{g}$. Next we define $$M: \mathfrak{g} \times \mathfrak{g} \rightarrow \C,$$ 
$$M(x,y) = m(\exp(x) \exp(y)^{-1}) \phi(\exp(x)) \phi(\exp(y)).$$ 
Observe that we have 
$$M(x,y) = m(\exp(x) \exp(y)^{-1}) \quad \mbox{for all} \quad x,y \in \mathrm{V} := \exp^{-1}(\Sigma).$$
According to \cite[Theorem 1.19]{LdlS}, we deduce
\begin{eqnarray*}
\big\| T_m: L_p(\V) \to L_p(\V) \big\|_{\mathrm{cb}} \!\! & \lesssim & \!\! \big\| S_m: S_p(\Sigma) \to S_p(\Sigma) \big\|_{\mathrm{cb}} \\ \!\! & = & \!\! \big\| S_M: S_p(\mathrm{V}) \to S_p(\mathrm{V}) \big\|_{\mathrm{cb}} \\ \!\! & \le & \!\! \big\| S_M: S_p(\R^n) \to S_p(\R^n) \big\|_{\mathrm{cb}}
\end{eqnarray*}
since $\mathfrak{g}$ is $n$-dimensional. Then, we get 
\begin{eqnarray} \label{Eq-EuclideanBound}
\lefteqn{\hskip-10pt \big\| T_m: L_p(\V) \to L_p(\V) \big\|_{\mathrm{cb}}} \\ \nonumber \hskip10pt \!\! & \le & \!\! C_p(\Omega) \sum_{|\gamma| \le [\frac{n}{2}]+1} \Big\| |x-y|^{|\gamma|} \Big\{ \big| \partial_x^\gamma M (x,y) \big| + \big| \partial_y^\gamma M (x,y) \big| \Big\} \Big\|_\infty.
\end{eqnarray}
This follows from Theorem \ref{Thm-HMS}. Then, we need to relate the Euclidean derivatives of $M$ with the Lie derivatives of the symbol $m: \G \to \C$ in a small neighborhood of the identity. To that end, now we claim that there exist smooth functions $a_{\alpha,\gamma}, b_{\alpha,\gamma} \in C_c^\infty(\R^n \times \R^n)$ supported by $\mathrm{U} \times \mathrm{U}$ for $|\alpha| \le |\gamma| \le [\frac{n}{2}]+1$ such that the following identities hold for every $|\gamma| \le [\frac{n}{2}]+1$
\begin{eqnarray} \label{Eq-EuclideanvsLie1}
\partial_x^{\gamma} M(x,y) \!\! & = & \!\!\!\! \sum_{|\alpha| \le |\gamma|} a_{\alpha, \gamma}(x,y) \, d^{\alpha} m \big( \text{exp}(x) \text{exp}(y)^{-1} \big), \\ \label{Eq-EuclideanvsLie2} \partial_y^{\gamma} M(x,y) \!\! & = & \!\!\!\! \sum_{|\alpha| \le |\gamma|} b_{\alpha, \gamma}(x,y) \hskip1pt \, d^{\alpha} m \big( \text{exp}(x) \text{exp}(y)^{-1} \big).
\end{eqnarray}
Before proving our claim, we show it implies Theorem A. Indeed, the Riemannian and Euclidean geometries are locally equivalent. In particular, since $M$ is supported in $\mathrm{U} \times \mathrm{U}$, we deduce that $L_R(\text{exp}(x) \text{exp}(y)^{-1})^{-|\alpha|} \sim |x - y|^{-|\alpha|} \lesssim |x - y|^{-|\gamma|}$ for every $x,y \in \mathrm{U}$ and every $|\alpha| \leq |\gamma|$. Therefore, combining this estimate with \eqref{Eq-EuclideanBound} and \eqref{Eq-EuclideanvsLie1} + \eqref{Eq-EuclideanvsLie2} we get the desired result. 

Let us now prove the claim. Given the definition of $M$, it is clearly true for $|\gamma| = 0$. We now justify it for $|\gamma| = 1$. By symmetry, we may just consider the $\mathrm{X}_1$-directional derivative in the variable $x$  
\begin{eqnarray*}
\partial_x^{e_1} M(x,y) \!\! & = & \!\! \lim_{s \rightarrow 0} \frac{M(x + s \mathrm{X}_1, y) - M(x,y)}{s} \\ \!\! & = & \!\! \frac{d}{ds}\Big|_{s=0}  m \big( \exp(x + s \mathrm{X}_1) \exp(y)^{-1} \big) \phi(\exp(x)) \phi(\exp(y)) \\ \!\! & + & \!\! m \big( \exp(x) \exp(y)^{-1} \big) \frac{d}{ds}\Big|_{s=0} \phi(\exp(x + s \mathrm{X}_1)) \phi(\exp(y)) \ =: \ \mathrm{A} + \mathrm{B}.
\end{eqnarray*}
The term $\mathrm{B}$ is already of the expected form in \eqref{Eq-EuclideanvsLie1} for $\alpha = 0$. Now observe that
$$\exp(x + s \mathrm{X}_1) \exp(y)^{-1} = \exp(x) \exp(y)^{-1} \exp(\mathrm{Z}_s),$$
where $\mathrm{Z}_s$ is given by
\begin{eqnarray*}
\mathrm{Z}_s(x,y) \!\! & = & \!\! \log \Big( \exp(y) \exp(x)^{-1} \exp(x + s \mathrm{X}_1) \exp(x)^{-1} \exp(x) \exp(y)^{-1} \Big) \\ \!\! & = & \!\! \exp(y) \exp(x)^{-1} \log \Big( \exp(x + s \mathrm{X}_1) \exp(x)^{-1} \Big) \exp(x) \exp(y)^{-1}.
\end{eqnarray*}
This is well defined for $s$ small enough since $\mathrm{supp} (M) \subset \mathrm{U} \times \mathrm{U}$. Hence
$$\frac{d}{ds}\Big|_{s=0}  m \big( \exp(x + s \mathrm{X}_1) \exp(y)^{-1} \big) \, = \, \lim_{s \rightarrow 0} \frac{\Phi_s(1) - \Phi_s(0)}{s} \, = \, \lim_{s \rightarrow 0} \frac{\Phi_s'(r(s))}{s}$$
for $\Phi_s(r) = m \big( \exp(x) \exp(y)^{-1} \exp(r \mathrm{Z}_s) \big)$ and some $0 < r(s) < 1$. Moreover 
\begin{eqnarray*}
\Phi_s'(r(s)) \!\! & = & \!\! \frac{d}{du}\Big|_{u = 0} m \big( \exp(x) \exp(y)^{-1} \exp(r(s) \mathrm{Z}_s) \exp(u \mathrm{Z}_s) \big) \\ \!\! & = & \!\! \partial_{\mathrm{Z}_s} m \big( \exp(x) \exp(y)^{-1} \exp(r(s) \mathrm{Z}_s) \big), 
\end{eqnarray*}
where $\partial_{\mathrm{Z}_s}$ denotes the left-invariant Lie derivative in the direction of $\mathrm{Z}_s$. Thus 
$$\mathrm{A} = \lim_{s \to 0} \Big\langle \frac{\mathrm{Z}_s(x,y)}{s}, \nabla m \big( \exp(x) \exp(y)^{-1} \exp(r(s) \mathrm{Z}_s) \big) \Big\rangle \phi(\exp(x)) \phi(\exp(y))$$ with $\nabla$ standing for the Lie gradient. We claim that 
\begin{equation} \label{Eq-Cinfty}
\lim_{s \to 0} \frac{\mathrm{Z}_s(x,y)}{s} \in C^\infty(\mathrm{U} \times \mathrm{U}).
\end{equation}
Therefore, since $\mathrm{Z}_s(x,y) \to 0$ as $s \to 0$ and $0 < r(s) < 1$, this implies that $\mathrm{A}$ is again of the expected form in \eqref{Eq-EuclideanvsLie1} with $|\alpha|=1$. In particular, \eqref{Eq-Cinfty} implies our claim for $|\gamma| = 1$. Higher order derivatives are dealt with in exactly the same way. In conclusion, all what is left to complete the proof of Theorem A is to justify claim \eqref{Eq-Cinfty}. Recalling the definition of $\mathrm{Z}_s(x,y)$ we are interested in computing $\Psi_s(x) = \log (\exp(x + s \mathrm{X}_1) \exp(x)^{-1})$. By the Hausdorff-Baker-Campbell formula we have
\begin{equation} \label{Eq-HBC}
\Psi_s(x) = \sum_{n \geq 1} \frac{(-1)^n}{n} \sum_{\begin{subarray}{c} r_j + t_j > 0 \\ (1 \le j \le n) \end{subarray}} \frac{\big\{ (x + s \mathrm{X}_1)^{r_1} (-x)^{t_1} \cdots (x+s\mathrm{X}_1)^{r_n} (-x)^{t_n} \big\}}{\sum_{j = 1}^n (r_j + t_j) \prod_{j = 1}^n r_j! t_j!},
\end{equation} 
where we use iterated Lie brackets 
$$\big\{ \mathrm{a}^{r_1} \mathrm{b}^{t_1} \cdots \mathrm{a}^{r_n} \mathrm{b}^{t_n} \big\} = \big[ \mathrm{d}_1, \big[ \mathrm{d}_2,  \big[ \ldots \big[\mathrm{d}_{m-1}, \mathrm{d}_m \big] \ldots \big] \big] \big]$$ with $m = \sum_{j \le n} r_j + t_j$ and where the first $r_1$ terms $\mathrm{d}_j$ equal $\mathrm{a}$, the next $t_1$ terms equal $\mathrm{b}$, the next $r_2$ terms equal $\mathrm{a}$ and so on. Next, note $[x + s \mathrm{X}_1,-x] = s[x,\mathrm{X}_1]$ and $$\{ (x + s \mathrm{X}_1)^{r_1} (-x)^{t_1} \cdots (x+s\mathrm{X}_1)^{r_n} (-x)^{t_n} \} = \pm s \hskip1pt [x,[x,[\ldots,\mathrm{X}_1] \ldots]]] + O(s^2),$$ unless $n=1$ and $r_1+t_1=1$. In this case we may find the additional terms $x$ for $(r_1,t_1)=(1,0)$ and $-x$ for $(r_1,t_1) = (0,1)$. However, these later terms have vanishing sum in \eqref{Eq-HBC} and we may ignore them. Thus, since $\Psi_s$ is absolutely convergent in a small neighborhood of $0$
$$\Psi_s(x) = O(s^2) + s \underbrace{\sum_{n \geq 1} \frac{(-1)^n}{n} \sum_{\begin{subarray}{c} r_j + t_j > 0 \\ (1 \le j \le n) \end{subarray}} a_{r_1,t_1,...,r_n,t_n} [x,[x,[\ldots,\mathrm{X}_1] \ldots]]]}_{\Pi(x)}$$ for some smooth function $\Pi$ near the origin. Therefore
\begin{eqnarray*}
\lim_{s \to 0} \frac{\mathrm{Z}_s(x,y)}{s} \!\! & = & \!\! \exp(y) \exp(x)^{-1} \lim_{s \to 0} \frac{\Psi_s(x)}{s} \exp(x) \exp(y)^{-1} \\ 
\!\! & = & \!\! \exp(y) \exp(x)^{-1} \Pi(x) \exp(x) \exp(y)^{-1} \in \ C^\infty(\mathrm{U} \times \mathrm{U}).
\end{eqnarray*}
This justifies our claim \eqref{Eq-Cinfty}, which in turn completes the proof of Theorem A. \fin

\section{\bf Simple Lie groups} \label{Sect-Simple}

Consider a simple Lie group $\G$. Define $\tau_\G = \mathrm{d}_\G / [\frac12 (\dim \G + 1)]$, where $\mathrm{d}_\G$ was introduced in \cite{Mau} as a crucial parameter in the volume growth of $\mathrm{Ad}$-balls up to a logarithmic factor 
\begin{equation} \label{Eq-Mau}
\mu \big\{ g \in \G: \|\mathrm{Ad}_g\| \le R \big\} \, \sim_{\log R} R^{\mathrm{d}_\G}.
\end{equation}
If $\mathcal{C}$ is the convex hull of the weights of the adjoint representation, Maucourant proved in \cite[Theorem 1.1]{Mau} that $\mathrm{d}_\G$ is the unique positive real number such that the sum of positive roots with multiplicities divided by $\mathrm{d}_\G$ lies on the boundary of $\mathcal{C}$. As $\mathcal{C}$ is a cone and the weights of the adjoint representation are precisely the roots, we get $$\sum_{\alpha \in \Sigma_+} \dim (\mathfrak{g}_\alpha) \alpha \in \Big( \sum_{\alpha \in \Sigma_+} \dim (\mathfrak{g}_\alpha) \Big) \mathcal{C}.$$ In particular, since $\mathfrak{g} = \bigoplus_{\alpha \in \Sigma \cup \{0\}} \mathfrak{g}_\alpha$ and $\Sigma = \Sigma_+ \cup \Sigma_-$, we obtain
\begin{equation} \label{Eq-Mik1}
\mathrm{d}_\G \, \le \, \sum_{\alpha \in \Sigma_+} \dim (\mathfrak{g}_\alpha) = \frac{\dim \G - \dim \mathfrak{g}_0}{2} \quad \mbox{and} \quad \tau_\G \, \le \, \frac{\dim \G}{2[\frac12 (\dim \G + 1)]} \, \le \, 1.
\end{equation}
Next, let $\Phi: \R_+ \to \R_+$ be a smooth function satisfying $$\Phi(x) \approx \begin{cases} x & \mbox{when} \ x < 1, \\ x^{\tau_\G} & \mbox{when} \ x \to \infty. \end{cases}$$ Then we set $$L_\G(g) := \Phi \big( \|\mathrm{Ad}_g - \mathrm{Ad}_e\| \big).$$ Around the identity, $L_\G$ behaves like the weight $g \mapsto \|\mathrm{Ad}_g - \mathrm{Ad}_e\|$ which is known to be locally Euclidean. On the other hand, asymptotically we obtain the behavior $$L_\G(g) \approx \|\mathrm{Ad}_g\|^{\tau_\G}.$$ This choice was already justified in \cite{PRS} for $\G = S \hskip-1pt L_n(\R)$ and $\tau_\G = \frac12$. 

\begin{lemma} \label{Lem-Exponential}
Given $\phi \in \mathcal{C}^1(\G \setminus \{e\})$ and $\beta>1/\tau_\G$
$$\sup_{\|\mathrm{X}\|=1} L_\G(g)^{\beta} \big| \partial_\mathrm{X} \phi (g) \big| \, \le \, 1 \ \Rightarrow \ L_\G(g)^\beta |\phi(g) - \alpha| \le C_\beta \quad \mbox{for some $\alpha \in \C$}.$$ Here we write $\partial_\mathrm{X}$ to denote the left-invariant Lie derivative in the direction $\mathrm{X} \in \mathfrak{g}$.
\end{lemma}

\dem We claim that $\phi - \alpha \in \mathcal{C}_0(\G)$ for some $\alpha \in \C$. This claim gives the statement. Indeed, according to the KAK decomposition, every $g \in \G$ factorizes as $g = k_1 k_2 k_2^{-1} \exp(\mathrm{Z}) k_2 = k \exp(s \mathrm{X})$ for some vector $\mathrm{Z}$ in the Cartan algebra, some unit vector $\mathrm{X} \in \mathfrak{g}$ and some $k \in \mathrm{K}$. By assumption, we obtain
\begin{eqnarray*}
| \phi(g) - \alpha| \!\!\!\! & = & \!\!\!\! \Big| \sum_{k \ge 1} \phi \big( g \exp((k-1)\mathrm{X}) \big) - \phi \big( g \exp(k \mathrm{X}) \big) \Big| \\ \!\!\!\! & = & \!\!\!\! \Big| \sum_{k \ge 1} \partial_\mathrm{X} \phi \big( g \exp(s_k \mathrm{X}) \big) \Big| \le \sum_{k \ge 1} L_\G \big(k \exp((s+s_k) \mathrm{X}) \big)^{-\beta}
\end{eqnarray*}
for some $s_k \in (k-1,k)$. Moreover, we have $\mathrm{Ad}_{\exp \mathrm{Z}} = \exp (\mathrm{ad}_\mathrm{Z})$ and $$\|\mathrm{Ad}_g\| = \exp \|\mathrm{Z}\| \quad \mbox{for} \quad g = k_1 \exp (\mathrm{Z}) \, k_2 \quad \mbox{and} \quad \|\mathrm{Z}\| = \max_{\alpha \in \Sigma} \alpha(\mathrm{Z}).$$ In particular, K-biinvariance of the map $g \mapsto \| \mathrm{Ad}_g \|^{\tau_\G}$ and $s \mathrm{X} = \mathrm{Ad}_{k_2}(\mathrm{Z})$ give
$$| \phi(g) - \alpha| \le \sum_{k \ge 1} L_\G \Big( \exp \big( \frac{s+s_k}{s} \mathrm{Z} \big) \Big)^{-\beta} \approx \sum_{k \ge 1} e^{-\beta \tau_\G (s+s_k)} \lesssim L_\G(g)^{-\beta}.$$
Let us now justify the claim. Using once more $g = k_1 \exp(\mathrm{Z}) k_2$ and the surjectivity of the exponential map onto $\mathrm{K}$, we get that $k_j = \exp (\mathrm{A}_j)$ for some $\mathrm{A}_j \in \mathfrak{k}$ with $\|\mathrm{A}_j\| \le 2\pi$. Under this factorization, we have $L_\G(g) = \exp (\tau_\G \|\mathrm{Z}\|)$ and we conclude
\begin{eqnarray} \label{eq-A2}
\lefteqn{\hskip-20pt \Big| \phi \big( \exp(\mathrm{A}_1) \exp(\mathrm{Z}) \exp(\mathrm{A}_2) \big) - \phi \big( \exp (\mathrm{A}_1) \exp(\mathrm{Z}) \big) \Big|} \\ \nonumber
\hskip50pt & = & \|\mathrm{A}_2\| \, \Big| \partial_{\hskip-2pt \frac{\mathrm{A}_2}{\|\mathrm{A}_2\|}} \phi \big( \exp(\mathrm{A}_1) \exp(\mathrm{Z}) \exp (r \mathrm{A}_2) \big) \Big| \, \le \, 2\pi \exp(-\beta \tau_\G \|\mathrm{Z}\|)
\end{eqnarray}
for some $0 < r < 1$. Similarly, let us note that $\exp(\mathrm{A}_1) \exp(\mathrm{Z}) = \exp(\mathrm{Z}) \mathrm{w}$ for $\mathrm{w} = \exp(- \mathrm{Z}) \exp(\mathrm{A}_1) \exp(\mathrm{Z}) = \exp (\mathrm{Y})$ where $\mathrm{Y} = \exp(- \mathrm{Z}) \mathrm{A}_1 \exp(\mathrm{Z})$ belongs to $\mathfrak{g}$. Therefore, the following identity holds for some $r \in (0,1)$ 
\begin{equation} \label{eq-A1}
\Big| \phi \big( \exp(\mathrm{A}_1) \exp(\mathrm{Z}) \big) - \phi \big(\exp(\mathrm{Z}) \big) \Big| = \|\mathrm{Y}\| \, \Big| \partial_{\hskip-2pt \frac{\mathrm{Y}}{\|\mathrm{Y}\|}} \phi \big( \exp(\mathrm{Z}) \exp (r \mathrm{Y}) \big) \Big|.
\end{equation}
Since $\|\mathrm{Y}\| = \|\mathrm{Ad}_{\exp \mathrm{Z}}(\mathrm{A}_1)\| = \|\exp(\mathrm{ad}_\mathrm{Z}(\mathrm{A}_1))\| \le 2\pi \exp \|\mathrm{Z}\|$ and $$L_\G(\exp(\mathrm{Z}) \exp (r \mathrm{Y})) = L_\G(\exp(r \mathrm{A}_1) \exp(\mathrm{Z})) = \exp (\tau_\G \|\mathrm{Z}\|),$$ the above quantity is bounded by $2\pi \exp(-(\beta\tau_\G-1)\|\mathrm{Z}\|)$, which decreases to $0$ for any $\beta > 1/\tau_\G$ as $\mathrm{Z} \to \infty$. According to \eqref{eq-A2} and \eqref{eq-A1}, it suffices to prove that $\phi - \alpha \in \mathcal{C}_0$ when restricted to elements $g=\exp \mathrm{Z}$ in the abelian part of the KAK decomposition. To prove it, consider the Euclidean function $\rho(\mathrm{Z}) = \phi(\exp \mathrm{Z})$ and fix $\mathrm{U}$ in the Cartan algebra. Then   
\begin{eqnarray*}
\langle \nabla \rho (\mathrm{Z}), \mathrm{U} \rangle & = & \lim_{s \to 0} \frac{\rho(\mathrm{Z} + s\mathrm{U}) - \rho(\mathrm{Z})}{s} \\ & = & \lim_{s \to 0} \frac{\phi(\exp(\mathrm{Z}) \exp(s \mathrm{U})) - \phi(\exp(\mathrm{Z}))}{s} \, = \,  \partial_\mathrm{U} \phi(\exp(\mathrm{Z})).
\end{eqnarray*} 
In particular, the Euclidean function $\rho$ satisfies
$$\sup_{\|\mathrm{U}\| = 1}  \big| \langle \nabla \rho (\mathrm{Z}), \mathrm{U} \rangle \big| \le \sup_{\|\mathrm{X}\| = 1}  \big| \partial_\mathrm{X} \phi (\exp(\mathrm{Z})) \big| \le L_\G(\exp(\mathrm{Z}))^{-\beta} \le \exp(-\beta \tau_\G \|\mathrm{Z}\|).$$ This readily implies that $\rho$ has a limit $\alpha$ at infinity and the same holds for $\phi$. \fin

\demAAA Since $L_\G$ is locally Euclidean around the identity, we know from Theorem A that the statement holds for Fourier symbols $m$ with small enough support. Of course, by a simple cut and paste argument, the same holds for symbols $m$ supported in a compact set $\Sigma$, as long as we admit that the cb-norm of $T_m$ could depend on $\Sigma$ and $p$. To justify the general statement, let $\Gamma$ be a cocompact lattice in $\G$ with fundamental domain $\Delta$ containing the identity $e$ in its interior. Consider a relatively compact neighborhood of the identity $\Omega$ containing the closure of $\Delta$ and satisfying that $e$ does not lie in the interior of $\gamma \Omega$ for any $\gamma \in \Gamma \setminus \{e\}$. Moreover, we may clearly impose that there exists an upper bound $\mathrm{M}$ for the number of translates $\gamma \Omega$ that may overlap with a fixed one $\gamma_0 \Omega$. Given $\phi \in \mathcal{C}_c^{\infty}(\G)_+$ supported in $\Omega$ and identically $1$ over $\Delta$, define $$\Phi_\gamma(g) = \frac{\phi(\gamma g)^2}{\sum_{\rho \in \Gamma} \phi(\rho g)^2} \quad \mbox{for each} \quad \gamma \in \Gamma.$$ The $\Phi_\gamma$'s form a smooth partition of unity in $\G$. Decompose $$m = \summ_{\gamma}\Phi_\gamma^{\frac12} m_\gamma \quad \mbox{where we take} \quad m_\gamma = \Phi_\gamma^{\frac12} m.$$ Note that the Fourier multipliers associated to the symbols $\sqrt{\Phi_\gamma}$ are completely bounded in $L_1(\V)$ (with the same cb-norm) and completely contractive in $L_2(\V)$. Indeed, it is clear by construction that $\sqrt{\Phi_\gamma} \in \mathcal{C}_c^{\infty}(\G)_+$ and or claim on the $L_1$-bound |equivalently the $L_\infty$-bound| follows by transference. Namely by smooth-cutoff our symbol into finitely many pieces, we may assume that our smooth symbol has small support around the identity (translation-invariance). By Bo\.zejko-Fendler transference from Remark \ref{Rem-Linfty-FSTransf}, we may consider the corresponding Herz-Schur multiplier, or even the Schur multiplier lifted to the Lie algebra via the exponential map as we did in the proof of Theorem A. It is at this point where we use that the support of our symbol is smooth enough. Once we have lifted the multiplier to the Lie algebra, our claim follows from \cite[Lemma 2.12]{PRS}. Next consider the linear map $$\Lambda: \ell_p(\Gamma; L_p(\V)) \to L_p(\V) \quad \mbox{with} \quad \Lambda \big( (f_\gamma)_{\gamma \in \Gamma} \big) = \summ_\gamma T_{\sqrt{\Phi_{\gamma}}}(f_\gamma).$$ Its cb-boundedness for $p=1$ follows from the triangle inequality and the (uniform) $L_1$ cb-boundedness of the multipliers in the above sum. The cb-contractivity for $p=2$ follows from Plancherel theorem and the finite overlapping property imposed on the translates $\gamma \Omega$ by construction. In particular, by complex interpolation $\Lambda$ is cb-bounded as well for $1 < p < 2$. Since $T_m(f) = \Lambda((T_{m_\gamma}(f))_\gamma)$, we get
$$\big\| T_m: L_p(\V) \to L_p(\V) \big\|_{\mathrm{cb}} \le \Big( \sum_{\gamma \in \Gamma} \big\| T_{m_\gamma}: L_p(\V) \to L_p(\V) \big\|_{\mathrm{cb}}^p \Big)^{\frac1p}.$$
\vskip-8pt \noindent Conjugating with the translation by $\gamma$, we may replace $m_\gamma$ by its left translate $M_\gamma(g) = m(\gamma^{-1} g) \sqrt{\Phi_e(g)}$. Then, the local behavior above yields for $\sigma = [\frac{n}{2}]+1$
\begin{eqnarray} \label{Eq-ellpBound}
\hskip20pt \|T_m\|_{\mathrm{cb}(L_p)}^p \!\!\!\! & \le & \!\!\!\! C_p \sum_{|\beta| \le \sigma} \sum_{\gamma \in \Gamma} \, \sup_{g \in \Omega} L_\G(g)^{p|\beta|} \big| d_g^\beta M_\gamma(g) \big|^p \\ \nonumber
\!\!\!\! & \le & \!\!\!\! C_p \sum_{|\beta| \le \sigma} \Big( \sup_{g \in \Omega} L_\G(g)^{p|\beta|} \big| d_g^\beta M_e(g) \big|^p + \sum_{\gamma \neq e} \, \sup_{g \in \Omega} \big| d_g^\beta M_\gamma(g) \big|^p \Big).
\end{eqnarray} 
\vskip-8pt \noindent Next, Leibniz rule and left-invariance of Lie differentiation give 
$$\big| d_g^\beta M_\gamma(g) \big| \, \le \, \sum_{\rho \le \beta} \big| d_g^{\rho} m(\gamma^{-1} g)  d_g^{\beta - \rho} \sqrt{\Phi_e} (g) \big| \, \lesssim \, \sum_{\rho \le \beta} \big| d_g^{\rho}m (\gamma^{-1} g) \big|.$$
According to the HM condition imposed in Theorem B, this shows that the $e$-term in the upper bound for \eqref{Eq-ellpBound} is controlled by the HM-norm of $m$. Moreover, since $e$ does not lie in the interior of $\gamma \Omega$ for any $\gamma \neq e$, a similar bound applies for the other $\gamma$-terms, but this is not enough to bound the whole sum. At this point we use Lemma \ref{Lem-Exponential}. In other words, according to it and the statement, every Lie derivative of $m$ with order $|\beta| \le \sigma$ decays asymptotically like $L_\G^{-\sigma}$. Here is the point where we use the assumption $n \ge 2/\tau_\G$, to ensure that $[\frac{n}{2}]+1 > 1/\tau_\G$. This implies
\begin{eqnarray*}
\sum_{\gamma \neq e} \sup_{g \in \Omega} \big| d_g^\beta M_\gamma(g) \big|^p \!\! & \lesssim & \!\! \sum_{\gamma \neq e} \, \sup_{g \in \Omega} \Big( \sum_{\rho \le \beta} \big| d_g^{\rho}m (\gamma^{-1} g) \big| \Big)^p \\ \!\! & \lesssim & \!\! \sum_{\gamma \neq e} \,  L_\G(\gamma)^{-\sigma p} \ \approx \ \sum_{\gamma \neq e} \,  \| \mathrm{Ad}_\gamma\|^{-\sigma \tau_\G p} \\ \!\! & = & \!\! \sum_{\gamma \neq e} \,  \| \mathrm{Ad}_\gamma\|^{- \mathrm{d}_\G p ([\frac{n}{2}]+1/[\frac{n+1}{2}])} \ \le \ \sum_{\gamma \neq e} \,  \| \mathrm{Ad}_\gamma\|^{- \mathrm{d}_\G p}.
\end{eqnarray*}
Next, we know from \eqref{Eq-Mau} that the volume of $\mathrm{Ad}$-balls of radius $R$ in $\G$ grows as $R^{\mathrm{d}_\G}$ up to a logarithmic factor. Thus, summing in dyadic $\mathrm{Ad}$-coronas we see that $\mathrm{d}_\G$ is the critical integrability index for $g \mapsto \|\mathrm{Ad}_g\|$ in $L_1(\G)$. Since $p>1$, the exponent $\mathrm{d}_\G p$ is above the critical index and the sum above is bounded by $C_p$. \fin

\begin{remark} \label{Rem-Optimal}
\emph{Note that $[\frac{n}{2}]+1 = [\frac{n+1}{2}]$ for odd $n=\dim \G$. In particular, the presence of $p>1$ is necessary in our argument. The proof in \cite{PRS} for $\G = \SL$ was simpler (only the triangle inequality is needed) since the HM condition there imposed derivatives up to order $[n^2/2]+1$ with $\dim \SL = n^2-1$. However, it seems that there is no more room for improvement in Theorem B. It is intriguing that this \lq\lq cut and paste" argument stops working precisely below the optimal HM condition since derivatives up to $[\dim \G/2]$ would give an exponent below the critical index. Is the asymptotic behavior in Theorem B optimal?}
\end{remark}

\begin{remark}
\emph{After Theorem B was made public, Martijn Caspers has proved an interesting result in \cite{Cas} which gives lower asymptotic decay rates of certain class of ($\mathrm{K}$-biinvariant) smooth symbols, which fail to admit a singular behavior around the identity. This is certainly interesting information in line with Theorem B. The statement resembles Calder\'on-Torchinsky interpolation theorem \cite{CT}.}
\end{remark}

\section{\bf HMS multipliers in groups} \label{Sect-Groups}

In this section, we investigate an extension of Theorem \ref{Thm-HMS} for Schur multipliers on locally compact groups with H\"ormander-Mikhlin-Schur (HMS) conditions in terms of group cocycles, we also consider weaker Sobolev regularity. Then we apply it to generalize the Fourier multiplier results from \cite{JMP1}.

\subsection{Singular multipliers in groups} 

Let $\G$ be a locally compact group equipped with a left Haar measure $\mu$ and a $n$-dimensional cocycle $\beta$. More precisely, there exists an action $\alpha: \G \to G \hskip-1pt L_n(\R)$ which does not need to be orthogonal, for which $\beta: \G \to \R^n$ satisfies the cocycle law 
\begin{equation} \label{Eq-CocycleLaw}
\alpha_g(\beta(h)) = \beta(gh) - \beta(g).
\end{equation}
Consider a smooth function $\phi: \R^n \to \R_+$ supported by $\mathrm{B}(0,2)$ and identically $1$ over $\mathrm{B}(0,1)$. Set $\psi(x) = \phi(x) - \phi(2 x)$, which is supported by $\mathrm{B}(0,2) \setminus \mathrm{B}(0,\frac12)$ and thus $\sum_{j \in \Z} \psi(2^j\xi)$ equals $1$ for all $\xi \neq 0$. Let $M: \G \times \G \to \C$ and assume there exist $\beta$-lifts $$M_{r\beta}: \R^n \! \times \G \to \C \quad and \quad M_{c\beta}: \G \times \R^n \to \C$$ satisfying 
\begin{equation} \label{Eq-BetaLifts}
M(g,h) = M_{r\beta}\big( \beta(g^{-1}) - \beta(h^{-1}),h \big) = M_{c\beta} \big(g,\beta(h^{-1}) - \beta(g^{-1}) \big).
\end{equation}
Define its \textbf{Sobolev HMS$_{\beta \sigma}$-norm} of order $\sigma = n/2+\varepsilon$ as follows $$\hskip32pt \big\bracevert \hskip-2pt M \hskip-2pt \big\bracevert_{\hskip-3pt \beta \sigma} \hskip-1pt := \hskip-1pt  \inf_{\begin{subarray}{c} \beta-\mathrm{lifts} \\ M_{r\beta}, M_{c\beta}\end{subarray}} \sup_{\begin{subarray}{c} j \in \Z \\ g,h \in \G \end{subarray}} \big\| \psi M_{r\beta}(2^{-j}\cdot,h) \big\|_{W_{2\sigma}} + \big\| \psi M_{c\beta} (g,2^{-j}\cdot) \big\|_{W_{2\sigma}}.$$ Here $W_{2\sigma}$ denotes the inhomogeneous $L_2$-Sobolev space in $\R^n$ of order $\sigma$. The above infimum will be assumed to be infinity when there are no $\beta$-lifts $M_{r\beta}$ or $M_{c\beta}$.

\begin{theorem}[General HMS multipliers] \label{ThmHMS}
\hskip-3pt Let $\G$ be a locally compact group equipped with a nonnecessarily orthogonal $n$-dimensional cocycle $\beta: \G \to \R^n$. Consider $M: \G \times \G \to \C$, $1 < p < \infty$ and $\sigma = n/2 + \varepsilon$ for some $\varepsilon > 0$. Then, the following inequality holds $$\big\| S_M \hskip-2pt: S_p(\G) \to S_p(\G) \big\|_{\mathrm{cb}} \, \le \, C_{\varepsilon} \frac{p^2}{p-1} \big\bracevert \hskip-2pt M \hskip-2pt \big\bracevert_{\hskip-2pt \beta \sigma}.$$
\end{theorem} 

\ske The proof is similar to \cite[Theorem A']{CGPT} and we will just highlight here the main differences. We shall divide the argument into four blocks:

\noindent \textbf{1) Reduction to BMO.} Let $\RR_\G$ be the von Neumann algebra of matrix-valued functions $L_\infty(\R^n) \bar\otimes \mathcal{B}(L_2(\G))$ and consider $\pi_\beta: \mathcal{B}(L_2(\G)) \to L_\infty (\RR_\G)$ to be the $*$-homomorphism  
$$\pi_\beta (A)(x) := \Big( \exp \big( 2\pi i \big\langle x, \beta(h^{-1}) - \beta(g^{-1}) \big\rangle \big) A_{gh} \Big).$$ 
Define $\mathrm{BMO}_\beta$ on $\B(L_2(\G))$ using $\pi_\beta$ as done in \cite{CGPT} via
$$\|A\|_{\mathrm{BMO}_{\beta}^c} = \sup_{Q \in \Q} \Big\| \Big( \mean_Q \big| \pi_\beta(A)(x) - \pi_\beta(A)_Q \big|^2 dx \Big)^{\frac12} \Big\|_{\mathcal{B}(L_2(\G))} = \big\|\pi_\beta(A) \big\|_{\mathrm{BMO}_{\mathcal{R}_\G}^c},$$  where $\Q$ is the set of Euclidean balls and $f_Q$ is the $Q$-mean of $f$. The seminorm in $\mathrm{BMO}_{\beta}$ vanishes on matrices $A$ satisfying that $\pi_\beta(A)$ is constant. In other words, on $\Sigma_p(\G) = \big\{ A \in S_p(\G): A_{gh} = 0 \ \mbox{ for } \ \beta(g^{-1}) \neq \beta(h^{-1}) \big\}$. Let $E_p: S_p(\G) \to \Sigma_p(\G)$ be the projection, which is a Schur multiplier $(g,h) \mapsto \delta_0 \big( \beta(h^{-1}) - \beta(g^{-1}) \big)$. Note that \eqref{Eq-CocycleLaw} gives $$\delta_0 \big( \beta(h^{-1}) - \beta(g^{-1}) \big) = \delta_0 \big( \alpha_g(\beta(h^{-1}) - \beta(g^{-1})) \big) = \delta_0 ( \beta(gh^{-1}) ).$$ Hence, $E_p$ is a Herz-Schur multiplier with symbol $(g,h) \mapsto \delta_0(\beta(gh^{-1}))$. By Remark \ref{Rem-GralFSTransf}, it is a cb-contraction. Indeed, $g \mapsto \delta_0(\beta(g))$ defines a cb-contractive Fourier multiplier on $\V$. This follows since the subgroup $\G_\beta = \{g \in \G: \beta(g)=0\}$ is either open |therefore inducing a conditional expectation| or with vanishing Haar measure. In particular, it turns out that $S_p^\circ (\G) := (\mathrm{id} - E_p)(S_p(\G))$ is a cb-complemented subspace of $S_p(\G)$. Then, according to \cite{JM} or an straightforward adaptation of \cite[Theorem 1.1]{CGPT} we get
\begin{equation} \label{Eq-Interp}
\big[ \mathrm{BMO}_{\beta}, S_2^\circ(\G) \big]_{\frac2p} \simeq_{\mathrm{cb}} \, S_p^\circ(\G) \quad \mbox{with} \quad c_p \approx p.
\end{equation}
Then, it suffices to prove that
\begin{equation} \label{Eq-BMOEndpoint2}
\big\| S_M: S_\infty(\G) \to \mathrm{BMO}_{\beta} \big\|_{\mathrm{cb}} \ \lesssim \ \big\bracevert \hskip-2pt M \hskip-2pt \big\bracevert_{\hskip-2pt \beta \sigma}.
\end{equation}
Indeed, using \eqref{Eq-Interp} and assuming \eqref{Eq-BMOEndpoint2}, we get 
\begin{eqnarray*}
\big\| S_M: S_p \to S_p \big\|_{\mathrm{cb}} \!\! & \lesssim & \!\! \big\| E_p S_M \big\|_{\mathrm{cb}(S_p)} + \big\| E_p^\perp S_M \big\|_{\mathrm{cb}(S_p)} \\ \!\! & \lesssim & \!\! \big\| E_p S_M: S_p \to S_p \big\|_{\mathrm{cb}} + \frac{p^2}{p-1} \big\bracevert \hskip-2pt M \hskip-2pt \big\bracevert_{\hskip-2pt L \beta \sigma}.
\end{eqnarray*}
However, by identity \eqref{Eq-BetaLifts} the symbol of $E_pS_M$ equals 
$$\delta_0 \big( \beta(h^{-1}) - \beta(g^{-1}) \big) M(g,h) = \delta_0(\beta(gh^{-1})) M_{c\beta}(g, 0).$$ Thus, $E_pS_M$ is the composition of a Schur multiplier (with symbol depending on $g$ but not on $h$) with $E_p$. As $E_p$ is a complete contraction, the cb-norm of $E_pS_M$ is bounded by the $L_\infty$-norm of the map $g \mapsto M_{c\beta}(g,0)$, which is bounded above by  the $\mathrm{HMS}_{\beta \sigma}$-norm. Hence, it remains to prove the BMO endpoint inequality \eqref{Eq-BMOEndpoint2}.

\noindent \textbf{2) Twisted Fourier multipliers.} 
Define the maps on $L_2(\RR_\G)$
\begin{eqnarray*}
\widetilde{T}_{M_{r\beta}}(f)(x) & := & \Big( \int_{\R^n} M_{r\beta}(\xi,h) \widehat{f}_{gh}(\xi) e^{2\pi i \langle x, \xi \rangle}\, d\xi \Big), \\
\widetilde{T}_{M_{c\beta}}(f)(x) & := & \Big( \int_{\R^n} M_{c\beta}(g,\xi) \widehat{f}_{gh}(\xi) e^{2\pi i \langle x, \xi \rangle}\, d\xi \hskip1pt \Big).
\end{eqnarray*}
The intertwining identities $\pi_\beta^* S_M = \widetilde{T}_{M_{r\beta}} \pi_\beta^*$ and $\pi_\beta S_M = \widetilde{T}_{M_{c\beta}} \pi_\beta$ give 
\begin{eqnarray} 
\label{Eq-Transf3}
\big\| S_M: S_\infty(\G) \to \mathrm{BMO}_{\beta}^r \big\|_{\mathrm{cb}} \!\! & \le & \!\! \big\| \widetilde{T}_{M_{r\beta}}: L_\infty(\RR_\G) \to \mathrm{BMO}_{\RR_\G}^r \big\|_{\mathrm{cb}}, \\
\label{Eq-Transf4}
\big\| S_M: S_\infty(\G) \to \mathrm{BMO}_{\beta}^c \big\|_{\mathrm{cb}} \!\! & \le & \!\! \big\| \widetilde{T}_{M_{c\beta}}: L_\infty(\RR_\G) \to \mathrm{BMO}_{\RR_\G}^c \big\|_{\mathrm{cb}}.
\end{eqnarray}

\noindent \textbf{3) Calder\'on-Zygmund kernels.} Now we claim that
\begin{eqnarray} 
\label{Eq-HM4}
\hskip20pt \big\| \widetilde{T}_{M_{r\beta}}: L_\infty(\RR_\G) \to \mathrm{BMO}_{L\RR_\G}^r \big\|_{\mathrm{cb}} \!\! & \lesssim & \!\! \sup_{(j,h) \in \Z \times \G} \ \big\| \psi M_{r\beta}(2^j \cdot,h) \big\|_{W_{2\sigma}}, \\
\label{Eq-HM3}
\hskip20pt \big\| \widetilde{T}_{M_{c\beta}}: L_\infty(\RR_\G) \to \mathrm{BMO}_{L\RR_\G}^c \big\|_{\mathrm{cb}} \!\! & \lesssim & \!\! \sup_{(j,g) \in \Z \times \G} \ \big\| \psi M_{c\beta}(g,2^j \cdot) \big\|_{W_{2\sigma}}.
\end{eqnarray}
To prove these inequalities, we need a kernel representation for these maps
\begin{eqnarray*}
\widetilde{T}_{M_{c\beta}}(f)(x) \!\! & = & \!\! \Big( \int_{\R^n} M_{c\beta}(g,\xi) \widehat{f}_{gh}(\xi) e^{2\pi i \langle x, \xi \rangle}\, d\xi \Big) \\ [7pt]
\!\! & = & \!\! \Big( \int_{\R^n} \underbrace{\big[ M_{c\beta}(g,\cdot) \big]^{\vee}(x-y)}_{k_{c\beta}(g,x-y)} f_{gh}(y) \, dy \Big) \, = \, \int_{\R^n} K_{c\beta}(x-y) \cdot f(y) \, dy
\end{eqnarray*}    
with $K_{c\beta}: \R^n \setminus \{0\} \to \mathcal{B}(L_2(\G))$ the denoting \lq\lq diagonal-valued" function $$\big( (K_{c\beta}(x)h)(g) \big) = \big( k_{c\beta}(g,x) h(g) \big) = \lq\lq \mathrm{diag} \big( [ M_{c\beta}(g,\cdot) ]^{\vee}(x) \big) \cdot \mathrm{col} \big( h(g) \big)".$$
Here we abusively use the matrix-notation $K_{c\beta} = \mathrm{diag} \big( k_{c\beta}(g,\cdot) \big)$. The kernel $K_{c\beta}$ should be understood as an operator-valued distribution, although we shall only use the kernel representation above when $x \in \R^n \setminus \mathrm{supp}_{\R^n} f$, which is both meaningful and the right framework for noncommutative CZ theory \cite{CCP,JMP1,Pa1}. At this point the proof follows verbatim as in \cite[Theorem A']{CGPT}. \fin  

\begin{remark} \label{Rem-AnisotropicMikhlin}
\emph{The HMS$_{\beta \sigma}$ condition recovers \cite[Theorem A']{CGPT} for $\G = \R^n$ and $\beta$ the trivial cocycle. 
It holds as well under the (more demanding) Mikhlin conditions
\begin{equation*} 
\sup_{g,h \in \G} \sum_{|\gamma| \le [\frac{n}{2}] +1} \Big\| |\xi|^{|\gamma|} \Big\{ \big| \partial_\xi^\gamma M_{r\beta}(\xi,h) \big| + \big| \partial_\xi^\gamma M_{c\beta}(g,\xi) \big| \Big\} \Big\|_\infty < \infty.
\end{equation*}}
\end{remark}

\begin{remark} \label{Rem-Omega}
\emph{The only point in the proof of Theorem \ref{ThmHMS} which requires to use the cocycle law is to establish the boundedness of projection $E_p$ to show that the BMO endpoint inequality suffices. A quick look at the argument shows that Theorem \ref{ThmHMS} could have been formulated in many other measure spaces $(\Omega, \Sigma,\mu)$ other than groups, just providing a function $\beta: \Omega \to \R^n$ which yields bounded projections $E_p$.}
\end{remark}

\begin{remark} 
\emph{The $*$-homomorphism $\pi_\beta$ in the proof of Theorem \ref{ThmHMS} can be easily modified to produce $L_p$-analogues of \eqref{Eq-Transf3} + \eqref{Eq-Transf4} respectively. It just requires to introduce a sequence of Gaussians in the line of K. de Leeuw's theorem \cite{dL} as shown in \cite[Proposition 1.7]{PRS}. This is however unproductive, since the resulting twisted Fourier multipliers fail $L_p$-boundedness in general. This explains the necessity of working with both twists and illustrates a better behavior of (singular) Schur multipliers compared to the class of operator-valued Calder\'on-Zygmund operators.}
\end{remark}

\subsection{Herz-Schur multipliers} Fourier multipliers over general locally compact groups do require an additional insight to face the lack of a standard differential structure. A broader interpretation of tangent space was given in \cite{JMP1}. Namely, if $\beta: \G \to \R^n$ is a finite-dimensional orthogonal cocycle and $m: \G \to \C$ satisfies the lifting identities $m(g) = \widetilde{m} \circ \beta(g) = m' \circ \beta (g^{-1})$, the main discovery was that a H\"ormander-Mikhlin theory in group von Neumann algebras is possible in terms of the $\beta$-lifted symbols for unimodular groups
\begin{equation} \tag{HM$_\beta$} \label{Eq-HMbeta}
\big\| T_m\hskip-2pt: L_p(\V) \to L_p(\V) \big\|_{\mathrm{cb}} \lesssim C_p \hskip-4pt \sum_{|\gamma| \le [\frac{n}{2}] +1} \Big\| |\cdot |^{|\gamma|} \Big\{ \big| \partial_\xi^\gamma \widetilde{m} \big| + \big| \partial_\xi^\gamma m' \big| \Big\} \Big\|_\infty.
\end{equation}
Now we generalize \eqref{Eq-HMbeta} to nonorthogonal cocycles and nonunimodular groups.

\begin{corollary}[Herz-Schur multipliers] \label{Cor-Herz-Schur}
Let $M(g,h) = m(gh^{-1})$ be a Herz-Schur multiplier and assume $m(g) = \widetilde{m}(\beta(g)) = m'(\beta(g^{-1}))$ for some cocycle $\beta$ associated to a $($not necessarily orthogonal\hskip1pt$)$ action $\alpha: \G \curvearrowright \R^n$. Then, the following estimate holds for $1 < p < \infty$ 
\begin{eqnarray*}
\lefteqn{\big\| S_M: S_p(\G) \to S_p(\G) \big\|_{\mathrm{cb}}} \\ \!\! & \le & \!\! \frac{p^2}{p-1} \sup_{g \in \G} \sum_{|\gamma| \le [\frac{n}{2}] + 1} \Big\| |\xi|^{|\gamma|} \Big\{ \big| \partial_\xi^\gamma (\widetilde{m} \circ \alpha_g)(\xi) \big| + \big| \partial_\xi^\gamma (m' \circ \alpha_g)(\xi) \big| \Big\} \Big\|_\infty.
\end{eqnarray*}
Moreover, if the action $\alpha$ is orthogonal we may remove $\alpha_g$ from the above condition.
\end{corollary}

\dem Set $$M_{r\beta}(\xi,h) = m'(\alpha_h(\xi)) \quad \mbox{and} \quad M_{c\beta}(g,\xi) = \widetilde{m}(\alpha_g(\xi)).$$ By \eqref{Eq-CocycleLaw} we get $M_{r\beta} \big( \beta(g^{-1}) - \beta(h^{-1}),h \big) = M(g,h) = M_{c\beta} \big( g, \beta(h^{-1}) - \beta(g^{-1}) \big)$. Then, the statement follows from Theorem \ref{ThmHMS} for $\sigma=[\frac{n}{2}]+1$, see Remark \ref{Rem-AnisotropicMikhlin}. \fin

According to Fourier-Schur transference \eqref{Eq-FS}, Corollary \ref{Cor-Herz-Schur} recovers the Fourier multiplier theorem \eqref{Eq-HMbeta} for amenable groups. Moreover, by local transference Theorem \ref{Thm-LocalTransf}, the same holds for compactly supported multipliers on nonamenable groups. In addition, unlike in \cite{JMP1,JMP2} our approach gives similar results for Schur multipliers in nonunimodular groups. 

\begin{remark}
\emph{A Sobolev form of Corollary \ref{Cor-Herz-Schur} also holds via Theorem \ref{ThmHMS}.}
\end{remark}

\begin{remark} \emph{It is tempting to apply Corollary \ref{Cor-Herz-Schur} with the standard cocycle $\beta(g) = g-e$ in $S \hskip-1pt L_n(\R)$ associated to the nonorthogonal action $\alpha_g(A) = g \cdot A$ for any $n \times n$ matrix $A \in S_2^n$, as it was used in \cite{PRS}. However, the condition in Corollary \ref{Cor-Herz-Schur} becomes void in this case, since the supremum in $g \in S \hskip-1pt L_n(\R)$ turns out to be unbounded. This was expectable, since otherwise Theorem \ref{ThmHMS} would get in conflict with the rigidity results from \cite{LdlS,PRS}.}
\end{remark}

\section*{\bf Appendix -- Stratified Lie groups}

\renewcommand{\theequation}{A.1}
\addtocounter{equation}{-1}

The analysis on nilpotent Lie groups was born in the 70's with a special emphasis on harmonic analysis \cite{C, MrM, MuM, MRS, MS, SS} and hypoelliptic PDEs \cite{F77,FS74,RS76}. It gradually attracted more attention, first with Heisenberg groups and later with general stratified Lie groups, we refer to \cite{FR16, St93, Th98} for a historical overview. As explained above, the motivations for a similar research over group algebras come from geometric group theory and operator algebra, rather than PDEs. This (dual) setting primarily focuses on spectral multipliers, which arise by functional calculus on the subLaplacian. In our setting a spectral approach seems hopeless, since the subRiemannian metric is not conditionally negative. 

In this Appendix we shall discuss an open problem concerning a HM criterium for stratified Lie group algebras. The first step comes from an anisotropic form of Theorem \ref{ThmHMS}. Fix $L = (\ell_1, \ell_2, \ldots, \ell_n)$ with $\ell_j \ge 1$. Given $\lambda > 0$ consider the $L$-dilations and the anisotropic $\rho_L$-metric in $\R^n$ given by 
\begin{equation} \label{Eq-Anisotropic}
\delta_{L\lambda}(x) = \big( \lambda^{\ell_1} x_1, \lambda^{\ell_2} x_2, \ldots, \lambda^{\ell_n} x_n \big) \quad \mbox{and} \quad \rho_L(x) = \Big( \summ_j |x_j|^{\frac{2}{\ell_j}} \Big)^{\frac{1}{2}}.
\end{equation}
Note that $\rho_L(\delta_{L\lambda}(x)) = \lambda \rho_L(x)$. Given a multi-index $\gamma = (\gamma_1, \gamma_2, \ldots, \gamma_n)$ we set $$\{\gamma\} = \sum_{j=1}^n \ell_j \gamma_j.$$

\begin{Ctheorem}[Anisotropic HMS multipliers] \label{ThmHMSAni}
\hskip-3pt Let $1 < p < \infty$ and consider a Schur symbol $M: \R^n \times \R^n \to \C$. Then, the following inequality holds for any $n$-anisotropic metric $\rho_L$ $$\big\| S_M \big\|_{\mathrm{cb}(S_p(\R^n))} \, \lesssim \, C_{L,p} \sum_{|\gamma| \le [\frac{n}{2}] +1} \Big\| \rho_L(x-y)^{\{\gamma\}} \Big\{ \big| \partial_x^\gamma M(x,y) \big| + \big| \partial_y^\gamma M(x,y) \big| \Big\} \Big\|_\infty.$$
\end{Ctheorem} 

We shall omit the proof of the above result. It arises from a technical detour of the proof of Theorem \ref{ThmHMS}, which involves an anisotropic form of matrix-valued CZ theory. Its Toeplitz form is an anisotropic Euclidean HM theorem from \cite{BB}, which refines in turn a classical result in \cite{FR}. The presence of the powers $\{\gamma\}$ rather than $|\gamma|$ should be understood as an stratified HM condition. Namely, a derivative in the $j$-th direction is dealt with as a $\ell_j$-th order derivative in our HM condition. One is led to think that a similar statement should hold for stratified Lie groups. Instead of working with spectral multipliers |as in the dual literature alluded above| we consider Fourier multipliers with a stratified Mikhlin condition. More precisely, if $\G$ is an $n$-dimensional stratified Lie group and $(\ell_1, \ell_2, \ldots, \ell_n)$ are the homogeneous dilation weights rescaled so that $\min(\ell_j) = 1$, then $\gamma$-derivatives must be bounded by $\{\gamma\}$-powers of the subRiemannian metric.

\begin{Cconjecture}[HM for stratified Lie groups] \label{Conject} 
Let $\G$ be a $n$-dimensional stratified Lie group. Consider the subRiemannian length $L_\mathrm{SR}: \G \to \R_+$ associated to its homogeneous dilation with weights $(\ell_1, \ell_2, \ldots, \ell_n)$ as above. Then, the following inequality holds for any Fourier symbol $m: \G \to \C$ and $1 < p < \infty$ $$\big\| T_m: L_p(\V) \to L_p(\V) \big\|_{\mathrm{cb}} \, \lesssim \, \frac{p^2}{p-1} \sum_{|\gamma| \le [\frac{n}{2}]+1} \big\| L_\mathrm{SR}(g)^{\{\gamma\}} d_g^\gamma m(g) \big\|_\infty.$$
\end{Cconjecture} 

One is tempted to use Theorem \ref{ThmHMSAni} to prove Conjecture \ref{Conject} in a similar way as we proved the local Theorem A. Indeed, using the natural dilation maps on $\G$, the above conjecture easily reduces to prove it just for symbols supported in a small neighborhood of the identity. At this point, we may use Fourier-Schur transference and the exponential map as in the proof of Theorem A to lift the problem to the Lie algebra, where Theorem \ref{ThmHMSAni} could be applied. This anisotropic form of \eqref{Eq-EuclideanBound} certainly gives sufficient conditions for $L_p$-boundedness of $T_m$. However, the conjectured HM condition in terms of Lie derivatives requires to apply \eqref{Eq-EuclideanvsLie1} and \eqref{Eq-EuclideanvsLie2}. The problem comes from the fact that $|\alpha| \le |\gamma|$ does not imply $\{\alpha\} \le \{\gamma\}$ which is what is required in this case. 

A good illustration of the difficulties one can expect arises by a quick look at the 3D Heisenberg group $\mathrm{H}_3$. We shall represent it as the group of upper triangular $3 \times 3$ matrices with $1$'s in the diagonal. The homogeneous dilation weights are $(1,1,2)$. Using that $$\exp \left( \begin{array}{ccc} 0 & x_1 & x_3 \\ 0 & 0 & x_2 \\ 0 & 0 & 0 \end{array} \right) = \left( \begin{array}{ccc} 1 & \ x_1 & x_3 + \frac12 x_1 x_2 \\ 0 & \ 1 & x_2 \\ 0 & \ 0 & 1 \end{array} \right)$$ and following the proof of Theorem A, we get
\begin{eqnarray*} 
\partial_{x_1} M(x,y) \!\! & = & \!\!\!\! d^{e_1} m \big( \text{exp}(x) \text{exp}(-y) \big) + \Big(\frac12 x_2 - y_2 \Big) d^{e_3} m \big( \text{exp}(x) \text{exp}(-y) \big), \\ \partial_{x_2} M(x,y) \!\! & = & \!\!\!\! d^{e_1} m \big( \text{exp}(x) \text{exp}(-y) \big) - \Big(\frac12 x_1 - y_1 \Big) d^{e_3} m \big( \text{exp}(x) \text{exp}(-y) \big).
\end{eqnarray*}
This is the form of \eqref{Eq-EuclideanvsLie1} for first-order $x_1$ and $x_2$ derivatives. According to Theorem \ref{ThmHMSAni}, we expect $\partial_{x_1}M(x,y)$ to grow like $1/\rho_L(x-y)$, which is locally comparable to $1/L_{\mathrm{SR}}(\exp(x) \exp(-y))$. This estimate holds for the left-invariant Lie derivative $d^{e_1} m$ by the hypotheses of Conjecture \ref{Conject}, but $d^{e_3} m$ grows as $L_{\mathrm{SR}}^{-2}$ since $e_3$ lives in the second stratum. Similarly, we may not find the expected growth for $\partial_{x_2}M(x,y)$ to apply Theorem \ref{ThmHMSAni}. Only $\partial_{x_3}M(x,y)$ satisfies the expected growth, since we allow it to be quadratic in this case.

\vskip5pt

\noindent A.1. \textbf{In search of small coefficients.} The objection above could be easily solved if the coefficients of $d^{e_3}m$ were bounded by $\rho_L(x-y)$, but the functions $\frac12 x_1 - y_1$ and $\frac12 x_2 - y_2$ fail this estimate near $0$. Alternatively, we could use another local diffeomorphism $\phi: \R^3 \to \mathrm{H}_3$ $$\phi(x) = \left( \begin{array}{ccc} 1 & \ \phi_1(x) & \phi_3(x) \\ 0 & \ 1 & \phi_2(x) \\ 0 & \ 0 & 1 \end{array} \right)$$ other than the exponential map and the lifting $M_\phi(x,y) = m(\phi(x)\phi(y)^{-1})$ to see if this can be achieved. Using left-invariant Lie derivatives, it turns out that \eqref{Eq-EuclideanvsLie1} has the following form for first order derivatives 
$$\partial_{x_j} M_\phi(x,y) = \sum_{k=1}^3 a_{jk}(x,y) d^{e_j} m \big( \phi(x)\phi(y)^{-1} \big),$$ \vskip-27pt \null $$\Big( \underbrace{\partial_{x_j}\phi_1(x)}_{a_{j1}(x,y)}, \underbrace{\partial_{x_j}\phi_2(x)}_{a_{j2}(x,y)}, \underbrace{\partial_{x_j}\phi_3(x) + \big( \phi_1(y) - \phi_1(x) \big) \partial_{x_j} \phi_2(x) - \phi_2(y) \partial_{x_j}\phi_1(x)}_{a_{j3}(x,y)} \Big).$$ Note that we just need $|a_{j3}(x,y)| \le C \rho_L(x-y)$ for $j = 1,2$. Since $\phi_1$ is assumed to be Lipschitz and $|\xi| \le \rho_L(\xi)$ around $0$, we just need to work with a solution of the system of PDEs
\begin{eqnarray*}
\partial_{x_1}\phi_3(x) \!\! & = & \!\! \phi_2(x) \partial_{x_1}\phi_1(x), \\
\partial_{x_2}\phi_3(x) \!\! & = & \!\! \phi_2(x) \partial_{x_2}\phi_1(x).
\end{eqnarray*}
However, taking crossed derivatives we get $$\partial_{x_1}\phi_1(x) \partial_{x_2} \phi_2(x) = \partial_{x_1} \phi_2(x) \partial_{x_2}\phi_1(x).$$ This, together with the system above, readily implies that the determinant of the Jacobian $J_\phi$ is identically $0$. In other words, every solution of the system above fails to be locally invertible at every point $x$. Therefore, there is no local diffeomorphism $\phi: \R^3 \to \mathrm{H}_3$ satisfying the expected properties. A similar phenomenon occurs if we consider right-invariant Lie derivatives instead, or even if we mix them. 

\vskip5pt

\noindent A.2. \textbf{Perturbed differential operators.} Observe that $$\partial_{x_3}M(x,y) = d^{e_3}m \big( \exp(x) \exp(-y) \big).$$ Then, the following local perturbations of $(\partial_{x_1}, \partial_{x_2}, \partial_{x_3})$ $$\big( \Delta_{1y}, \Delta_{2y}, \Delta_{3} \big) = \Big( \partial_{x_1} + \frac12 y_2 \partial_{x_3}, \, \partial_{x_2} - \frac12 y_1 \partial_{x_3}, \, \partial_{x_3} \Big)$$ satisfy the identities 
\begin{eqnarray*}
\Delta_{1y}M(x,y) \!\!\! & = & \!\!\! d^{e_1}m \big( \exp(x) \exp(-y) \big) + \frac12 \big( x_2 - y_2 \big) d^{e_3}m \big( \exp(x) \exp(-y) \big), \\ \Delta_{2y}M(x,y) \!\!\! & = & \!\!\! d^{e_2}m \big( \exp(x) \exp(-y) \big) - \frac12 \big( x_1 - y_1 \big) d^{e_3}m \big( \exp(x) \exp(-y) \big). 
\end{eqnarray*}
According to the hypotheses in Conjecture \ref{Conject}, the right hand sides of both terms above grow as $L_{\mathrm{SR}}(\exp(x) \exp(-y))^{-1} \approx \rho_L(x-y)^{-1}$. Therefore, it suffices to prove a local form of Theorem \ref{ThmHMSAni} with the differential operators $\partial_{x_j}$ replaced by their local perturbations $\Delta_{jy}$. Similarly, we need to replace $\partial_{y_j}$ by local perturbations $\Delta_{jx}$ which are constructed in the same way from \eqref{Eq-EuclideanvsLie2} instead of \eqref{Eq-EuclideanvsLie1}. We can prove local perturbed forms of Theorem \ref{ThmHMS}, but it is unclear whether similar results hold for anisotropic metrics. Indeed, the operators $(\Delta_{1y}, \Delta_{2y}, \Delta_3)$ impose to work with different metrics for each $y$ in a local neighborhood of $0$. The difficulty comes from the fact that these metrics are nonequivalent when the underlying original metric is anisotropic. In conclusion, some additional ideas beyond \cite{CGPT} seem to be necessary to prove Conjecture \ref{Conject}. 

\begin{Cremark}
\emph{Some of our observations above apply as well to arbitrary stratified groups and essentially all of them are still valid for every stratified 2-step groups.}
\end{Cremark}

\noindent \textbf{Acknowledgement.} JP wants to express his gratitude to Mikael de la Salle for very valuable comments on a preliminary version (during a conference at Lorentz Center in 2021), which led to an improved presentation. We also want to thank \'Eric Ricard for relevant suggestions which have given rise to an improved presentation. The authors were partly supported by the Spanish Grant PID2019-107914GB-I00 \lq\lq Fronteras del An\'alisis Arm\'onico" (MCIN / PI J. Parcet) and Severo Ochoa Grant CEX2019-000904-S (ICMAT), funded by MCIN/AEI 10.13039/501100011033. J. Conde-Alonso was supported by the Madrid Government Program V PRICIT Grant SI1/PJI/2019-00514. E. Tablate was supported as well by Spanish Ministry of Universities with a FPU Grant with reference FPU19/00837. 


\enlargethispage{1cm}

\bibliographystyle{amsplain}


\enlargethispage{1cm}

\vskip4pt

\small

\noindent \textbf{Jos\'e M. Conde-Alonso} \hfill \textbf{Adri\'an M. Gonz\'alez-P\'erez} \\ 
Universidad Aut\'onoma de Madrid \hfill Universidad Aut\'onoma de Madrid
\\ Instituto de Ciencias Matem{\'a}ticas \hfill Instituto de Ciencias Matem{\'a}ticas
\\ \texttt{jose.conde@uam.es} \hfill \texttt{adrian.gonzalez@uam.es}

\vskip0pt

\noindent \textbf{Javier Parcet} \hfill \noindent \textbf{Eduardo Tablate} \\
Instituto de Ciencias Matem{\'a}ticas \hfill Instituto de Ciencias Matem{\'a}ticas 
\\ CSIC \hfill CSIC  
\\ \texttt{parcet@icmat.es} \hfill \texttt{eduardo.tablate@icmat.es}



\end{document}